\documentclass[11pt]{article}
\usepackage{amsmath}
\usepackage{amssymb}
\hsize \textwidth
\advance \hsize by -\marginparwidth
\evensidemargin \oddsidemargin
\advance\hoffset by 5mm

\newcommand{\no}{\nonumber}
\newcommand{\be}{\begin{equation}}
\newcommand{\ee}{\end{equation}}
\newcommand{\bi}{\begin{itemize}}
\newcommand{\ei}{\end{itemize}}
\newcommand{\br}{\begin{eqnarray}}
\newcommand{\er}{\end{eqnarray}}

\newcommand{\expE}{\mathbb{E}}
\newcommand{\qed}{$\Box$}

\newcommand{\norm}[1]{\lVert #1 \rVert}

\newcommand{\commentout}[1]{}

\newcommand{\Rm}{{\mathbb R}}
\newcommand{\Zm}{{\mathbb Z}}
\newcommand{\Em}{{\mathbb E}}
\newcommand{\Pm}{{\mathbb P}}
\newcommand{\var}{\mathop{\rm var}\nolimits}

\newtheorem{theorem}{Theorem}[section]
\newtheorem{prop}[theorem]{Proposition}

\newtheorem{lemma}[theorem]{Lemma}

\begin{document}

\title{A sublinear variance bound for
solutions of a random Hamilton-Jacobi equation}

\author{Ivan Matic\thanks{Mathematics Department, Duke University, Box 90320, Durham, North Carolina, 27708, USA. (matic@math.duke.edu). }\;\;\;  and \;\;James Nolen\thanks{Mathematics Department, Duke University, Box 90320, Durham, North Carolina, 27708, USA. (nolen@math.duke.edu). }}

\date{}
\maketitle


\begin{abstract}
We estimate the variance of the value function for a random optimal control problem. The value function is the solution $w^\epsilon$ of a Hamilton-Jacobi equation with random Hamiltonian  $H(p,x,\omega) = K(p) - V(x/\epsilon,\omega)$ in dimension $d \geq 2$. It is known that homogenization occurs as $\epsilon \to 0$, but little is known about the statistical fluctuations of $w^\epsilon$. Our main result shows that the variance of the solution $w^\epsilon$ is bounded by $O(\epsilon/|\log \epsilon|)$.  The proof relies on a modified Poincar\'e inequality of Talagrand.
\end{abstract}

\section{Introduction}

In this paper we study the random optimal control problem
\be
u(t,x,\omega) = \sup_{\gamma \in \mathcal{A}_{t,x}} \; g(\gamma(t)) -\mathcal{L}(\gamma,\omega), \quad x \in \Rm^d, \;\; t > 0 \label{unxdef}
\ee
in dimension $d \geq 2$, where the supremum is taken over the set of admissible paths
\[
\mathcal{A}_{t,x} = \{ \gamma \in W^{1,\infty}([0,t];\Rm^d) \;|\; \gamma(0) = x\}.
\]
The upper-semicontinuous payoff function $g:\Rm^d \to \Rm \cup \{-\infty\}$ is given. The cost functional $\mathcal{L}$ has the form
\[
\mathcal{L}(\gamma,\omega)= \int_0^t L(\gamma'(s) , \gamma(s),\omega)\,ds =  \int_0^t K(\gamma'(s)) + V(\gamma(s),\omega)\,ds,
\] 
where $K(p):\Rm^d \to [0,\infty)$ is convex and grows super-linearly in $|p|$. The function $V(x,\omega)$ is a scalar random field that is statistically stationary and ergodic with respect to certain translations in $x$. The parameter $\omega \in \Omega$ denotes a sample from a given probability space $(\Omega, \mathcal{F},\Pm)$. Thus, the value function $u(t,x,\omega)$ is random. Our main result shows that the variance of $u(x,t,\omega)$ grows only sublinearly in $t$ as $t \to \infty$.

Under certain conditions on $g$ and $L$, $u(t,x,\omega)$ is uniformly continuous and is a viscosity solution \cite{BCD} of the random Hamilton-Jacobi equation
\be 
\left\{\begin{array}{rcll}
&&u_t = H(Du,x,\omega), \quad x \in \Rm^d, \quad t > 0\vspace{3pt}\\
&&u(0,x)   = g(x), \quad x \in \Rm^d,\end{array}\right.
\ee
where $H(p,x) = K^*(p) - V(x)$, $K^*$ being the Legendre transform of $K$. 
For simplicity, consider the case where $g(x) = \eta \cdot x$ is a linear function. Then for each $\epsilon > 0$, the function $w^\epsilon(t,x,\omega) = \epsilon u(t/\epsilon,x/\epsilon,\omega)$ solves the initial value problem
\be 
\left\{\begin{array}{rcll}
&&w^\epsilon_t = H(Dw^\epsilon,\frac{x}{\epsilon},\omega), \quad x \in \Rm^d, \quad t > 0\vspace{3pt}\\
&&w^\epsilon(0,x)  =  \eta \cdot x = g(x) .\end{array}\right. \label{weps}
\ee
For certain Hamiltonians $H(p,x,\omega)$ which are convex in $p$, statistically stationary and ergodic with respect to translation in $x$, it is known \cite{RT,Soug1} that as $\epsilon \to 0$, homogenization occurs (see also \cite{AS, LS2} for alternative proofs and \cite{KV, KRV, LS1, Schw} for related results). This means that the functions $w^\epsilon(t,x,\omega)$ converge locally uniformly in $[0,\infty) \times \Rm^d$, as $\epsilon \to 0$ to the deterministic function $\bar w(t,x)$ which solves
\be 
\left\{\begin{array}{rcll}
&&\bar w_t = \bar H(D \bar w), \quad x \in \Rm^d, \quad t > 0\vspace{3pt}\\
&&\bar w(0,x)  =  g(x).\end{array}\right.
\ee
The function $\bar H(p):\Rm^d \to \Rm$ is called the effective Hamiltonian. We may think of this convergence as kind of law of large numbers for $w^\epsilon$, although the limit $\bar w$ and the effective Hamiltonian $\bar H$ are not determined by a simple averaging. Beyond this convergence result, relatively little is known about the properties of $\bar H$, about the rate of convergence $w^\epsilon \to \bar w$, or about the statistical behavior of $w^\epsilon - \expE[w^\epsilon]$, where $\expE[\cdot]$ denotes expectation with respect to the probability measure $\Pm$.  Our work pertains to this last issue: in terms of $w^\epsilon(t,x,\omega)$, our estimate on the variance of $u$ implies that $\var(w^\epsilon(t,x,\omega)) \leq C \epsilon/|\log \epsilon|$, as $\epsilon \to 0$.

Before stating our main result, let us make some definitions and assumptions more precise. We will suppose the random field $V(x,\omega)$ has the following special structure. Let $a<b$ be two real numbers. Let $\Omega = \{a,b\}^{\mathbb{Z}^d}$ be the set of all functions $\omega:\mathbb{Z}^d \to \{a,b\}$. Let the probability measure $\Pm$ be the shift-invariant product measure on $\Omega$ determined by $\Pm(\omega_k = a) = \alpha$ and $\Pm(\omega_k = b) = \beta$, for all $k \in \mathbb{Z}^d$, where $\alpha \in (0,1)$ and $\beta= 1-\alpha$. Thus the random variables $\{\omega_k\}_{k \in \mathbb{Z}^d}$ are independent and identically distributed. Now for $k \in \mathbb{Z}^d$, let $Q_k = k + [0,1)^d$ denote the unit cube with corner at the point $k$. Given $\omega \in \Omega$, define $V(x,\omega):\Rm^d \times \Omega \to \{a,b\}$ by
\be
V(x,\omega) = \sum_{k \in \mathbb{Z}^d} \omega_k \mathbb{I}_{Q_k}(x), \label{Vdef}
\ee
with $\mathbb{I}_{Q_k}$ is the indicator function for the set $Q_k$. Thus, $x \mapsto V$ is piecewise constant, taking values $a$ or $b$ on the unit cubes. By construction, the law of $V(x,\omega)$ is the same as that of $V(x + k,\omega)$ for any $k \in \mathbb{Z}^d$.  This precise construction of the field $V(x,\omega)$ is not essential for our result to hold. In particular, the function could be mollified so that it is uniformly continuous, or $V(x,\omega)$ could depend on the values of $\omega_k$ for $k$ in a bounded neighborhood of $x$. Nevertheless, the choice of $\Pm$ as the product measure on $\Omega = \{a,b\}^{\mathbb{Z}^d}$ is motivated by the main analytical tool presented below in Theorem \ref{talagrands_inequality}.

We suppose that $K:\Rm^d \to [0,\infty)$ is convex, $K(0) = 0$, and that 
\[
\lim_{|z| \to \infty} \frac{K(z)}{|z|} = +\infty.
\]
For the case of dimension $d = 2$ we will make use of an extra non-degeneracy condition: for some $\nu>1$, 
\be
K(z)\geq |z|^\nu \quad \forall \; z\in B_{1/2}(0). \label{nondeg2}
\ee
Given $V$ and $K$, let $L(p,x,\omega) = K(x) + V(x,\omega)$ and let $u$ be defined by (\ref{unxdef}). The following estimate of the variance of $u$ for large $t$ is our main result:

\begin{theorem}\label{fk_model_variance} Let $d \geq 2$. Let $x \in \Rm^d$ and suppose that $g:\Rm^d \to \Rm \cup \{-\infty\}$ satisfies
\be 
g(y) < g(x) + C_1(1 + |y - x|), \quad \forall y \in \Rm^d.
\ee
There is a constant $C > 0$, depending only on $C_1$, $K$, $\alpha$, $\beta$, and $|b - a|$, such that
\be
\var(u(t,x,\omega)) \leq C \frac{t}{\log t}, \quad \forall \; t \geq 2. \label{varbound1}
\ee
\end{theorem}

The main tool that we use to control the variance of $u(t,x,\omega)$, is the following theorem, which is a slight variation of an inequality of Talagrand (see \cite{talagrand_94}, Theorem 1.5). This result holds for product spaces of the form $\Omega_J = \{a,b\}^J$, $J$ being a finite set, and $\Pm$ being the product measure on $\Omega_J$ with marginals $\mathbb P(\omega_j = a)=\alpha \in (0,1)$ and $\mathbb P(\omega_j = b)=\beta = 1 - \alpha$, for all $j\in J$.  Let us define $\phi_j \omega$ to be the element of $\{a,b\}^J$ which is identical to $\omega$ except that the $j$-th component $\omega_j$ is opposite to $\omega_j$. That is, $\phi_j \omega = \omega'$, where $\omega_k' = \omega_k$ for $k \neq j$, and $\omega_j' \neq \omega_j$. For each random variable $f:\Omega_J \to \Rm$ define $\sigma_j f(\omega) = f(\phi_j \omega)$ and 
\[
\rho_jf(\omega)=\frac{\sigma_j f(\omega)-f(\omega)}{2}.
\]

\begin{theorem}\label{talagrands_inequality}
There is a constant $C>0$, independent of $|J|$, such that
\begin{eqnarray}\label{talagrand_ineq1}
\var(f)\leq C\sum_{j\in J}\frac{\|\rho_j f\|_2^2}{1+\log\frac{\|\rho_j f\|_2}{\|\rho_jf\|_1}}
\end{eqnarray}
holds for all $f \in L^2(\Omega_J)$.
\end{theorem}

The idea of using this inequality to estimate the variance of $f(\omega) = u(t,x,\omega)$ comes from the work of Benjamini, Kalai, and Schramm \cite{benjamini_kalai_schramm} who used this inequality to estimate the distance variance in first passage percolation, a problem which has some features similar to the control problem (\ref{unxdef}).  Specifically, they consider the length of minimal paths between two points in the integer lattice $\mathbb{Z}^d$ under a random metric. Each edge $e$ in the nearest-neighbor graph is assigned an independent random weight $\omega_e \in \{ a,b\}$, and the length of a path between two points $x,y \in \mathbb{Z}^d$ is defined as the sum of the edge weights along a path connecting $x$ and $y$.  They proved that $\var(d_\omega(0,v)) \leq C|v|/\log|v|$, where $d_\omega(0,v)$ is the length of the shortest path connecting $0$ and $v$. See \cite{BR} for some extensions of that result.  The main difficulty in applying the ideas of \cite{benjamini_kalai_schramm} to the present setting comes from the different structure of the cost functional $\mathcal{L}(\gamma,\omega)$, which necessitates more control on the optimizing paths.

As we have mentioned, for $d \geq 2$ there are relatively few results about the random fluctuation of $u(t,x,\omega)$ (as $t \to \infty$) or $w^\epsilon(t,x,\omega)$ (as $\epsilon \to 0$). In \cite{rezakhanlou_clt}, Rezakhanlou derived conditions under which a central limit theorem holds for $w^\epsilon(t,x)$ where $w^\epsilon$ is the solution of the Hamilton-Jacobi equation (\ref{weps}), i.e. whether $\epsilon^{-1/2}(w^\epsilon - \bar w)$ converges in law to some nontrivial stochastic process as $\epsilon \to 0$. In the case $d = 1$ those conditions can be verified for Hamiltonians having the form $H(p,x,\omega) = K(p) - V(x,\omega)$, and the limit distribution can be computed (see Corollary 2.6 in \cite{rezakhanlou_clt}). For $d \geq 2$, however, it is difficult to verify those conditions. Indeed, our result shows that we may have $\var(w^\epsilon) = o(\epsilon)$, which is less than what a CLT as in \cite{rezakhanlou_clt} would suggest. As this paper was being written, we learned of another work by Armstrong, Cardaliaguet, and Souganidis \cite{ACS}, who study the rate of convergence $w^\epsilon \to \bar w$. Our Theorem \ref{fk_model_variance} pertains to the variance of $w^\epsilon$, i.e the statistical error $w^\epsilon - \expE[w^\epsilon]$, but does not give an estimate of the bias $\expE[w^\epsilon] - \bar w$. 

The paper is organized as follows. In Section \ref{sec:optpaths} we derive some properties of the paths $\gamma$ which nearly optimize (\ref{unxdef}).  Section \ref{sec:maintheo} contains the main argument for the proof of Theorem \ref{fk_model_variance}. Section \ref{sec:tech} and Section \ref{sec:append} contain proofs so some technical estimates needed in Section \ref{sec:maintheo}.

\commentout{
The arguments used to proof Theorem \ref{fk_model_variance} extend to discrete time models. For example, consider the Frenkel-Kontorova type model
\begin{equation}
u(n,\omega) = \sup_{\gamma}  g(\gamma_n) - \mathcal{L}(\gamma,\omega), \label{unxdef_fk}
\end{equation}
where the supremum (\ref{unxdef_fk}) is taken over the set of admissible discrete sequences
\[
\mathcal{A}_{n} = \{\gamma = (\gamma_0, \gamma_1, \dots, \gamma_n) \;|\; \gamma_j \in \Rm^d,\;\;\gamma_0 = 0\}
\]
having fixed starting point $\gamma_0=0$, and the functional $\mathcal{L}$ is defined by
\begin{eqnarray}
\mathcal{L}(\gamma,\omega)= \sum_{k=0}^{n-1} K(\gamma_{k+1}-\gamma_k) + V(\gamma_k,\omega).
\end{eqnarray} 
}

\vspace{0.2in}

\noindent{\bf Acknowledgement} The work of JN was partially funded by grant DMS-1007572 from the US National Science Foundation.

\section{Properties of optimizing paths} \label{sec:optpaths}

Without loss of generality, let us suppose $x = 0$ and simply write $u(t,\omega) = u(t,0,\omega)$ and 
\[
\mathcal{A}_{t} = \mathcal{A}_{t,0} = \{ \gamma \in W^{1,\infty}([0,t];\Rm^d) \;|\; \gamma(0) = 0\}.
\]
So, we are studying the quantity
\be
u(t,\omega) = \sup_{\gamma \in \mathcal{A}_t} \; g(\gamma(t)) -\mathcal{L}(\gamma,\omega).
\ee
for some function $g:\Rm^d \to \Rm \cup \{-\infty\}$, with $g(0) \in \Rm$. For $\delta > 0$ and $\omega \in \Omega$, let $M_\delta(t,\omega)$ be the set of all paths $\gamma \in \mathcal{A}_t$ such that
\[
 g(\gamma(t)) - \mathcal{L}(\gamma,\omega)  \geq  u(t,\omega) - \delta.
\]
For each $\delta > 0$, this set is non-empty, and we refer to these paths as $\delta$-approximate optimizers. Observe that $M_{\delta_2}(t,\omega) \subset M_{\delta_1}(t,\omega)$ if $0 < \delta_2 < \delta_1$. If an optimal path $\gamma$ exists, meaning that $u(t,\omega) = g(\gamma(t)) - \mathcal{L}(\gamma,\omega)$, then it is certainly an approximate optimizer for any $\delta > 0$. In this section we derive some useful properties of approximate optimizers which will be used in the proof of Theorem \ref{fk_model_variance}.

\subsection*{Deterministic Bounds}

First, we have a few estimates which do not involve the random structure of the control problem.

\begin{lemma} \label{lem:upperC}
Let $\gamma \in M_\delta(t,\omega)$. Then for any $r_1,r_2 \in [0,t]$,
\begin{eqnarray}
\int_{r_1}^{r_2} L(\gamma'(s),\gamma(s),\omega)\,ds \leq (r_2 - r_1) b + (r_2 - r_1) K\left( \frac{\gamma(r_2) - \gamma(r_1)}{r_2 - r_1}\right) + \delta, \label{upperLsumC}
\end{eqnarray}
and
\begin{eqnarray}
\int_{r_1}^{r_2} L(\gamma'(s),\gamma(s),\omega)\,ds \geq (r_2 - r_1) a + (r_2 - r_1) K\left( \frac{\gamma(r_2) - \gamma(r_1)}{r_2 - r_1}\right). \label{lowerLsumC}
\end{eqnarray}
\end{lemma}

\vspace{0.2in}

\noindent{\bf Proof of Lemma \ref{lem:upperC}:} Given $\gamma \in M_\delta(t,\omega)$, define a new path $\hat \gamma  \in \mathcal{A}_{t}$ according to 
\[
\hat \gamma(s) = \gamma(r_1) + (s - r_1)\frac{\gamma(r_2) - \gamma(r_1)}{r_2 - r_1}\;\; s \in [r_1,r_2]
\]
and $\hat \gamma(s) = \gamma(s)$ for $s \notin [r_1,r_2]$. Thus we have replaced a section of $\gamma$ with a straight-line path connecting the same points. Since $\gamma \in \mathcal{M}_\delta(\omega)$, we must have $\mathcal{L}(\gamma,\omega) \leq \mathcal{L}(\hat \gamma,\omega) + \delta$. In particular,
\begin{eqnarray}
\int_{r_1}^{r_2} L(\gamma', \gamma, \omega)\,ds & \leq & \int_{r_1}^{r_2} L(\hat \gamma', \gamma, \omega)\,ds + \delta  \no \\
& \leq & (r_2 - r_1) b + (r_2 - r_1) K\left( \frac{\gamma(r_2) - \gamma(r_1)}{r_2 - r_1} \right).
\end{eqnarray}
This proves (\ref{upperLsumC}).  The lower bound (\ref{lowerLsumC}) follows from Jensen's inequality, the convexity of $K$, and the fact that $V(x,\omega) \geq a$.  \hfill \qed

\vspace{0.2in}

\begin{lemma} \label{lem:finspeed2C}
Suppose that $g:\Rm^d \to \Rm \cup \{-\infty\}$, $g(0) \in \Rm$, and
\be
g(y) < g(0) + C_1(1 + |y|), \quad \forall \;\; y \in \Rm^d. \label{finspeed2C}
\ee
There is a constant $R$ depending only on $K$, $C_1$, and $b - a$ such that
\[
|\gamma(t_2) - \gamma(t_1)| \leq R(1 + |t_1 - t_2|), \quad \forall \; t_2,t_1 \in [0,t]
\]
holds for all paths $\gamma \in M_1(t,\omega)$ and all $t >1$.
\end{lemma}

\vspace{0.2in}

\noindent{\bf Proof of Lemma \ref{lem:finspeed2C}:}  We first show there is a constant $R_0$ depending only on $K$, $C_1$, and $b - a$ such that
\be
|\gamma(t) - \gamma(0)| \leq tR_0 \label{finspeed}
\ee
holds for all $\gamma \in  M_1(t,\omega)$ and all $t \geq 1$. Define the path $\hat \gamma(s) = \gamma(0) = 0$ for all $s \in [0,t]$. We have
\begin{eqnarray}
u(t,\omega) \geq g(0) - \mathcal{L}(\hat \gamma,\omega) \geq g(0) - t b .
\end{eqnarray}
By (\ref{lowerLsumC}), we also have the lower bound
\[
\mathcal{L}(\gamma,\omega)   \geq t K\left(\frac{\gamma(t) - \gamma(0)}{t}\right) + a t.
\]
Since $\gamma \in M_1(t,\omega)$, we may combine these two estimates with $u(t,\omega) \leq 1 + g(\gamma(t)) - \mathcal{L}(\gamma,\omega)$ to conclude
\[
K\left(\frac{\gamma(t) - \gamma(0)}{t}\right) \leq  \frac{1}{t} + (b - a) + \frac{g(\gamma(t)) - g(0)}{t} \leq  1 + (b - a) + C_1\left( \frac{1 + |\gamma(t)|}{t} \right).
\]
Since $K(p)$ grows super-linearly in $|p|$, (\ref{finspeed}) follows.

\vspace{0.1in}

Next, consider $\gamma$ at integer times $k \in [1,t-1] \cap \mathbb{Z}$.  We will show that there is a constant $R_1$, independent of $t > 1$, such that at least one time $k \in [1,t-1] \cap \mathbb{Z}$ must satisfy both
\begin{eqnarray}  \label{good_terms}
|\gamma(k) - \gamma(k-1)|&\leq& R_1 \;\;\mbox{ and }\;\;  |\gamma(k+1) - \gamma(k)|
\leq R_1. 
\end{eqnarray}
Arguing by way of contradiction, let us suppose (\ref{good_terms}) does not hold. Then $|\gamma(j+1) - \gamma(j)| > R_1$ must hold for at least $t/3$ of the times $j \in [1,t-1] \cap \mathbb{Z}$.  This implies that
\[
u(t,\omega) \leq 1+ g(\gamma(t)) - \mathcal{L}(\gamma,\omega) \leq 1+ g(\gamma(t)) - at - \frac{t}{3} \min_{|q| \geq R_1} K(q).
\]
On the other hand, by Lemma \ref{lem:upperC} and (\ref{finspeed}), we know that
\[
u(t,\omega) \geq  g(\gamma(t)) - bt - t K(\frac{\gamma(t) - \gamma(0)}{t}) - 1 \geq   g(\gamma(t)) - bt -t \max_{|q| \leq R_0} K(q)  - 1
\]
holds for all $\gamma \in M_1(t,\omega)$. Combining these two bounds we obtain
\[
\frac{1}{3} \min_{|q| \geq R_1} K(q) \leq 1 +  (b- a) +  \max_{|q| \leq R_0} K(q).
\]
If $R_1 > R_0$ is sufficiently large (depending only on $b-a$, $R_0$, and $K$) this forces a contradiction.  So, (\ref{good_terms}) must hold.

\vspace{0.1in}

Now we conclude the proof.  Let $R_2 > R_1$, and suppose that for some $t_1,t_2 \in [0,t]$ with $1 \leq |t_2 - t_1| \leq 2$ we have $|\gamma(t_2) - \gamma(t_1)|\geq R_2$. Let $k \in [1,t-1] \cap \mathbb{Z}$ be such that (\ref{good_terms}) holds. Without loss of generality, we may suppose $k + 1 \leq t_1 < t_2$. Consider the path $\hat \gamma$ defined by
\[
\hat \gamma(s)=
\left\{\begin{array}{cl}\gamma(s),& \mbox{ for }s \in [0,k-1] \cup [t_2,t],\\
       \gamma(k-1) + (s - k + 1) (\gamma(k+1) - \gamma(k-1)),&\mbox{ for }s\in [k-1,k], \no \\
\gamma(s + 1), &\mbox{ for }s\in [k,t_1-1], \no
\end{array}\right.
\]
and for $s \in [t_1 - 1,t_2]$
\[
\hat \gamma(s) = \gamma(t_1) + (\gamma(t_2) - \gamma(t_1))\frac{s - t_1 + 1}{t_2 - t_1 + 1}.
\]
Then we have 
\br  
\mathcal L(\hat \gamma)-\mathcal L( \gamma)&\leq & 4(b - a) + \int_{k-1}^{k} K( \hat \gamma'(s)) \,ds + \int_{t_1-1}^{t_2} K( \hat \gamma'(s)) \,ds  \no \\
&  & - \int_{k-1}^{k+1} K(\gamma'(s)) \,ds - \int_{t_1}^{t_2} K(\gamma'(s)) \,ds \no \\
& \leq & 4( b - a )+ K( \gamma(k+1) - \gamma(k-1)) +  (t_2 - t_1 + 1) K( \frac{\gamma(t_2) - \gamma(t_1)}{t_2 - t_1 + 1}) \no \\
& & -  2K\left(\frac{\gamma(k+ 1) - \gamma(k-1)}{2}\right) - (t_2 - t_1)K\left(\frac{\gamma(t_2) - \gamma(t_1)}{t_2 - t_1}\right) \no \\
&\leq&  M +  (t_2 - t_1 + 1) K\left( \frac{\gamma(t_2) - \gamma(t_1)}{t_2 - t_1 + 1}\right)  - (t_2 - t_1)K\left(\frac{\gamma(t_2) - \gamma(t_1)}{t_2 - t_1}\right), \label{Lbig}
\end{eqnarray}
where 
\[
M= 4(b-a)+\max_{|z|\leq 2 R_1} K(z).
\]
Let $\Delta t = t_2 - t_1$ and $\sigma = (\Delta t + 1)/(\Delta t)$ and $z = (\gamma(t_2) - \gamma(t_1))/(t_2 - t_1 + 1)$. The inequality (\ref{Lbig}) has the form
\be
\mathcal L(\hat \gamma)-\mathcal L( \gamma) \leq M  + (\Delta t + 1)K(z) - \Delta t K(\sigma z). \label{MR2bound}
\ee
The properties of $K$ (convexity and super-linear growth) imply that if $R_2$ sufficiently large, then
\[
\inf_{|z| \geq R_2/3} K(\sigma z) - \sigma K(z) > M + 1.
\]
Applying this at (\ref{MR2bound}) we conclude $\mathcal L(\gamma')-\mathcal L( \gamma) < 1$, which contradicts the fact that $\gamma \in M_1(t,\omega)$. Therefore, we must have $|\gamma(t_2) - \gamma(t_1)| \leq R_2$ if $1 \leq |t_1 - t_2| \leq 2$. This and the triangle inequality now imply the desired result for all $t_1,t_2 \in [0,t]$. \hfill\qed

\vspace{0.2in}

\subsection*{Important cubes}
Our method of estimating the variance of $u$ involves bounding the random variable $|\sigma_j u - u|$. So, we must understand when changing the value of $\omega_j$ leads to a large change in the value of $u(t,\omega)$. Given a path $\gamma \in \mathcal{A}_{t}$ and an index $j \in \mathbb Z^d$, define 
\[
\pi_j(\gamma) = |\{ s \in [0,t] \;|\; \gamma(s) \in Q_j \}|,
\]
which is the total time that the path $\gamma$ occupies the cube $Q_j$. Observe that for any path $\gamma \in \mathcal{A}_t$, we have
\begin{equation}
\mathcal{L}(\gamma,\phi_j \omega) \leq \mathcal{L}(\gamma, \omega) + (b - \omega_j) \pi_j(\gamma). \label{sigdif1}
\end{equation}
In particular, if $\gamma \in M_\delta(t,\omega)$, then
\begin{eqnarray}
\sigma_j u(t,\omega) =  u(t,\phi_j \omega) & \geq & g(\gamma(t)) - \mathcal{L}(\gamma,\phi_j \omega) \no \\
&  \geq & g(\gamma(t)) - \mathcal{L}(\gamma, \omega) - (b - \omega_j) \pi_j(\gamma)  \no \\
& \geq & u(t,\omega) - (b - \omega_j) \pi_j(\gamma) - \delta. \label{sigdif}
\end{eqnarray}
From this we deduce that if $\omega_j = b$ or if there is $\gamma \in M_\delta(t,\omega)$ for which $\pi_j(\gamma) = 0$, then it must be true that $u(t,\omega) - \sigma_j u(t,\omega) \leq  \delta$. On the other hand, this also shows that if $u(t,\omega) - \sigma_j u(t,\omega) > \delta$, then $\omega_j = a$ and $\pi_j(\gamma) > 0$ must hold for all $\gamma \in M_\delta(t,\omega)$. This motivates the following definition. We say that the cube $Q_j$ is {\em important} if $\omega_j = a$ and for some $\delta > 0$ we have
\be
\pi_j(\gamma) > 0, \;\; \forall\;\gamma \in M_{\delta}(t,\omega). \label{impdef}
\ee
Observe that if (\ref{impdef}) holds for some $\delta > 0$, then it also holds for all $\delta' \in (0,\delta]$. So, $Q_j$ is important if $\omega_j = a$ and for $\delta$ sufficiently small every $\delta$-approximate optimizer spends time in cube $Q_j$. Let ${\mathcal{I}}_j \subset \Omega$ denote the event that the cube $Q_j$ is important:
\br
{\mathcal{I}}_j & = & \left \{ \omega \in \Omega \;|\; \omega_j = a;\;\; \exists \delta > 0 \;\;\text{such that}\;\; \pi_j(\gamma) > 0 \;\; \forall\;\gamma \in M_\delta(t,\omega) \right\} \label{Ajdef} \\
& = & \bigcup_{n \geq 1}  \left \{ \omega \in \Omega \;|\; \omega_j = a;\;\;  \pi_j(\gamma) > 0 \;\; \forall\;\gamma \in M_{1/n}(t,\omega) \right\}.  \label{AjdefC}
\er
The above analysis shows that
\be
\{\omega \in \Omega \;|\; u(t,\omega) > \sigma_j u(t,\omega)\} \subset {\mathcal{I}}_j, \quad \forall\;j \in \Zm^d \label{PsigA1}
\ee
so we have
\begin{equation}
\mathbb P( \sigma_j u < u) \leq \Pm({\mathcal{I}}_j). \label{PsigA}
\end{equation}
Observe that $\Pm({\mathcal{I}}_j)$ depends on $t$, in addition to $j$.

It will be useful to further classify some cubes as {\it very important}. To this end, we define a set of cubes
\[
N(\delta,\omega) = \bigcup_{\gamma \in M_\delta(\omega)} \{ k \in \mathbb{Z}^d \;|\; \pi_k(\gamma) > 0 \}.
\]
This is the set of all cubes visited by some path $\gamma \in M_{\delta}(t,\omega)$. Next, we define the event ${\mathcal{I}}_j^+ \subset {\mathcal{I}}_j \subset \Omega$ that cube $Q_j$ is {\em very important}:
\br
{\mathcal{I}}^+_j = \{ \omega \in {\mathcal{I}}_j \;|\;  \exists \, \delta > 0 \;\; \text{such that} \;\; \omega_\ell = b\;\;\; \;\forall \;\ell \in N(\delta,\omega)\setminus \{j\} \}.
\er
On this event, $Q_j$ is an important cube, and for any other cube $Q_\ell$ visited by a path $\gamma \in M_\delta$, we have $\omega_\ell = b$, if $\delta$ is sufficiently small. On the event ${\mathcal{I}}^-_j = {\mathcal{I}}_j \setminus {\mathcal{I}}_j^+$, cube $Q_j$ is important but not very important: for any $\delta >0$ we can find a path $\gamma \in M_\delta(\omega)$ such that $\gamma$ passes through another cube $Q_\ell \neq Q_j$, on which $\omega_\ell = a$.  The following lemma shows that the only way for $(u - \sigma_j u)^21_{{\mathcal{I}}_j}$ to be large is if $Q_j$ is very important.

\vspace{0.2in}

\begin{lemma} \label{lem:notveryimportant}
There is a constant $C_0> 0$, depending only on $K$ and $|b-a|$, such that 
\[
\Pm\left(\{ \omega \;|\;\;(u - \sigma_j u)^2 1_{{\mathcal{I}}^-_j} < C_0 \}\right) = 1
\]
holds for all $t \geq 1$ and $j \in \Zm^d$.
\end{lemma}
\noindent {\bf Proof of Lemma \ref{lem:notveryimportant}:} 
If $\omega_j = b$, then $\omega \notin {\mathcal{I}}_j^-$, so obviously $(u - \sigma_j u)^21_{{\mathcal{I}}_j^-} = 0$.  Hence, we may assume $\omega_j = a$ and $\omega \in {\mathcal{I}}_j^-$. When $\omega_j = a$, we clearly have $\sigma_j u \leq u$, since $\mathcal{L}(\gamma,\phi_j \omega) \geq \mathcal{L}(\gamma,\omega)$ in this case. So, we must bound $u - \sigma_j u$ from above. 

Consider an approximate optimizer $\gamma \in M_\delta(\omega)$ for some $\delta \leq 1$.  If $\pi_j(\gamma) \leq C$ then $u- \sigma_j  u \leq (b - a)C + \delta$ according to (\ref{sigdif}). So, we must consider the possibility that $\pi_j(\gamma) > 0$ is large. Since $\omega \in {\mathcal{I}}_j^-$, we may assume the path  $\gamma$ also passes through another cube $Q_\ell$, with $\ell \neq j$, for which $\omega_\ell=a$. We will construct a new path $\hat \gamma$ such that $\pi_j(\hat \gamma) \leq 1$ and $\mathcal{L}(\hat \gamma,\omega) \leq \mathcal{L}(\gamma,\omega) + C$. The two paths $\hat \gamma$ and $\gamma$ will have the same starting and ending points. This implies that the difference 
$u - \sigma_j u$ is bounded by a constant, since by (\ref{sigdif1}) we have 
\begin{eqnarray}
\sigma_j u = u(t,\phi_j \omega) & \geq & g(\hat \gamma(t)) - \mathcal{L}(\hat \gamma, \phi_j \omega) \nonumber \\
 & \geq & g(\hat \gamma(t)) - \mathcal{L}(\hat \gamma, \omega) - (b - a) \pi_j(\hat \gamma) \nonumber \\
& \geq & g(\hat \gamma(t)) - \mathcal{L}(\gamma, \omega) - C - (b - a) \pi_j(\hat \gamma) \geq u(t,\omega) - C - (b-a) - \delta.
\end{eqnarray}
Suppose that $[t_1,t_2]$ is the smallest interval containing all $s$ for which $\gamma(s) \in Q_j$. We may assume $t_2 - t_1 \geq \pi_j(\gamma) > 1$.  Suppose that $\gamma(t_3) \in Q_\ell$ where 
$\ell \neq j$ and $\omega_\ell=a$. We may suppose that $t_3 > t_2$ (the case $t_3 < t_1$ is similar). Define the new path $\hat \gamma$ as follows:
\begin{enumerate}
\item[(i)] For $s \in [0,t_1]$, let $\hat \gamma(s) = \gamma(s)$.
\item[(ii)] For $s \in [t_1,t_1 + 1]$, let $\hat \gamma(s) = \gamma(t_1) + (s - t_1) (\gamma(t_2) - \gamma(t_1))$.
\item[(iii)] For $s \in [t_1 + 1,t_1 + 1 + (t_3 - t_2)]$, let $\hat \gamma(s) = \gamma(s - t_1 - 1 + t_2)$.
\item[(iv)] For $s \in [t_1 + 1 + (t_3 - t_2), t_3]$, let $\hat \gamma(s) = \gamma(t_3)$.
\item[(v)] For $s \in [t_3,t]$, let $\hat \gamma(s) = \gamma(s)$.
\end{enumerate}
Much of $\hat \gamma$ is just a linear reparameterization of $\gamma$, and we have
\br
&& \int_0^t L(\hat \gamma'(s), \hat \gamma(s),\omega)\,ds - \int_0^t L( \gamma'(s), \gamma(s),\omega)\,ds  \leq \int_{t_1}^{t_1 + 1} K(\hat \gamma'(s))\,ds \leq K(\gamma(t_2) - \gamma(t_1)). \no
\er
Since $|\gamma(t_2) - \gamma(t_1)|$ is bounded by the diameter of cube $Q_j$, we have $\mathcal{L}(\gamma',\omega) \leq \mathcal{L}(\gamma,\omega) + C$. \hfil \qed


\section{Proof of Theorem \ref{fk_model_variance}} \label{sec:maintheo}

In this section we prove Theorem \ref{fk_model_variance}. As we have mentioned, the main argument is similar to that of \cite{benjamini_kalai_schramm}. In particular, it is convenient to average $u(t,\omega)$ over a random shift of the environment.

\subsection*{Random shifting of the environment}
We now consider an augmented probability space $\tilde \Omega = \Omega \times \Omega_1$ with product measure $\tilde \Pm = \Pm \times P_1$, and we introduce a random function $h(\omega_1):\Omega_1 \to \mathbb{Z}^d$ to define a random shift of the environment. For $(\omega,\omega_1) \in \tilde \Omega$, let us define
\begin{eqnarray}
\tilde u(t,\omega,\omega_1) =  u(t,\tau_{h(\omega_1)} \omega) = \sup_{\gamma \in \mathcal{A}_{t}} g(\gamma(t) ) -  \mathcal{L}(\gamma + h(\omega_1),\omega) \label{ftilde}
\end{eqnarray}
where $\gamma + h(\omega_1)$ denotes the shifted path $t \mapsto \gamma(t) + h(\omega_1)$.  We define $M_\delta(\omega,\omega_1) = M_\delta(\tau_{h(\omega_1)} \omega)$ to be the set of paths $\gamma \in \mathcal{A}_t$ for which
\[
\tilde u(t,\omega,\omega_1) \leq g(\gamma(t)) - \mathcal{L}(\gamma + h(\omega_1),\omega) + \delta.
\]
We construct $\Omega_1$, $P_1$, and $h$ in such a way that $|\tilde u(t,\omega,\omega_1) - u(t,\omega)| = o(\sqrt{t})$, and for this reason an estimate of $\var (\tilde u)$ that is sublinear in $t$ will imply a sublinear bound for $\var(u)$.

The random shift $h(\omega_1)$ will lie in the set $[0,m)^d \subset \Rm^d$ where $m=\left\lfloor t^{\zeta}\right\rfloor$, for some positive $\zeta < 1/2$. For $d \geq 3$, it will suffice to choose $\zeta \in (\frac{1}{d}, \frac{1}{2})$. For $d = 2$, we will require that $\zeta \in \left(\frac{\nu-1}{2\nu-1},1/2\right)$, where $\nu$ was defined by the non-degeneracy condition (\ref{nondeg2}). Denote by $P_0$ the product probability measure on the set $\Omega_0 =\{a,b\}^{m^2}$ and having marginal distribution $P_0'(a) = \alpha$, $P_0'(b) = \beta = 1 - \alpha \in (0,1)$. The following statement is Lemma 3 from \cite{benjamini_kalai_schramm}, so we omit the proof:

\begin{lemma}\label{lemma3_benjamini_kalai_schramm} There exists a constant $C>0$ independent of $m= \left\lfloor t^{\zeta}\right\rfloor$, and a function $\tilde h:\Omega_0\to \{0,1,\dots, m-1\}$ for which the  following two conditions hold:
\begin{enumerate}
 \item[(i)] $P_0(\tilde h=i)\leq \frac C{m}$ for all 
$i\in [0,m-1]$, and
\item[(ii)] for every $x,y\in \Omega_0$ that differ in at 
most one coordinate, the difference between $\tilde h(x)$ and $\tilde h(y)$ satisfies $$|\tilde h(x)-\tilde h(y)|\leq 1.$$
\end{enumerate}
\end{lemma}

Define the set $\Theta = \Theta_t =\{1,2,\dots, d\}\times \{1,2,\dots, m^2\}
.$ Let $\Omega_1=\{a,b\}^{\Theta}$, and let $P_1$ be a uniform probability measure on $\Omega_1$.   Each $\omega_1\in\Omega_1$ can  be written as $\omega_1=(\omega_1^1,\omega_1^2,\dots, \omega_1^d)$, where each $\omega_1^i$ is a binary sequence of length $m^2$. Let $\vec e_i$ denote the $i$-th coordinate vector. Define
$h(\omega_1)=\sum_{i=1}^d  \tilde h(\omega_1^i)\vec e_i$. There exists a constant $C>0$ independent on $m$ such that for each $x\in\{0,1,\dots, m-1\}^d$ one has  $P_1 (h=x)\leq\frac C{m^{d}}$. Moreover, if $\omega_1$ and $\omega_1'$ differ in exactly one coordinate then we have $|h(\omega_1)-h(\omega_1')|\leq 1$. Given the space $\tilde \Omega = \Omega\times \Omega_1$ with the product measure $\tilde \Pm = \mathbb P\times P_1$ on $\Omega\times \Omega_1$ defined in this way, we now consider the function $\tilde u$ defined by (\ref{ftilde}).

\vspace{0.2in}

\begin{lemma} \label{lem:fftilde}
There is a constant $C > 0$ such that 
\begin{equation}
|u(t,\omega)-\tilde u(t,\omega,\omega_1)|\leq C|h(\omega_1)| \leq C t^{\zeta}, \quad \forall \; t > 1 \label{ftildefdif}
\end{equation}
holds $\tilde \Pm$-almost surely, and
\begin{equation}
\var u \leq C \var \tilde u + C t^{2 \zeta}, \quad \forall\;\; t > 1. \label{fftildevar}
\end{equation}
\end{lemma} 
\noindent{\bf Proof of Lemma \ref{lem:fftilde}:} We will prove that $|u(t,\omega)-\tilde u(t,\omega,\omega_1)|\leq C |h(\omega_1)|$, $\tilde \Pm$-almost surely,
for some constant $C>0$ independent of $m$ and $t$. Given a path $\gamma \in M_\delta(\omega)$, we can modify it to construct an approximate optimizer for $\tilde u(t,\omega,\omega_1)$, thus estimating $\tilde u(t,\omega, \omega_1) - u(t,\omega)$ from above. However, we cannot simply shift $\gamma$ by $-h(\omega_1)$, since we must preserve the starting and ending points. 

Suppose $|h(\omega_1)| \leq \sqrt{d} \kappa$, with $\kappa \in [1, m] \cap \mathbb{Z}$. Fixing a path $\gamma \in M_\delta(\omega)$, we define the new path $\hat \gamma$ in the following way: 
\begin{enumerate}
 \item[(i)] For $r \in [0, \kappa]$, set $\hat \gamma(r) = \gamma(0) + \left(\frac{r}{\kappa}\right)(\gamma(2\kappa) - h(\omega_1) - \gamma(0))$. 
 \item[(ii)] For $r \in [\kappa,t - \kappa]$, set $\hat \gamma(r) = \gamma(r + \kappa) - h(\omega_1)$. 
\item[(iii)] For $r\in [t-\kappa, t]$, set $\hat \gamma(r)= \gamma(t) + (r-t)\frac{h(\omega_1)}{\kappa}$.
\end{enumerate}

We now verify that the path $\hat \gamma$ yields the desired bound on $\tilde u(t,\omega,\omega_1)-u(t,\omega)$ and $\var u$. Since $\gamma \in M_\delta(\omega)$, we have
\begin{eqnarray*}
\tilde u(t,\omega,\omega_1)-u(t,\omega)&\geq &
- \int_{0}^{\kappa} L(\hat \gamma'(s),\hat \gamma(s) + h(\omega_1),\omega)\,ds -
\int_{t-\kappa}^{t} L(\hat \gamma'(s),\hat \gamma(s) + h(\omega_1),\omega) \,ds \no \\
&& + \int_0^{2\kappa} L(\gamma'(s), \gamma(s),\omega) \,ds - \delta\\
&\geq& - C \kappa \left( 1 + \sup_{|z|\leq 2 R}K(z) \right) - \delta,
\end{eqnarray*}
where $C$ is a positive real number that depends only on $K$, $|b-a|$. 
In a similar way we prove that 
$\tilde u(t,\omega,\omega_1)- u(t,\omega) \leq C\kappa + \delta$. Recalling that $m = \left\lfloor t^{\zeta}\right\rfloor$, we obtain (\ref{ftildefdif}). Therefore 
\begin{eqnarray*}
\var u&=& \tilde \expE[(u-\tilde \expE(u))^2]\\
&=&\tilde \expE\left[ \left(\tilde u-\tilde \expE(\tilde u)  + (u  - \tilde u) - \tilde \expE(u-\tilde u) \right)^2\right] \\
&\leq & 3 \var \tilde u+ 3\tilde \expE\left[\left((u-\tilde u)-\tilde \expE(u-\tilde u)\right)^2\right] \\
&\leq& 3\var\tilde u +12C^2m^2,
\end{eqnarray*}
which is (\ref{fftildevar}).

 \hfill $\Box$
 
\vspace{0.2in}

\subsection*{Variance estimate for $\tilde u$}

Given Lemma \ref{lem:fftilde} and the choice of $\zeta < 1/2$, we now wish to establish a bound of order $t/\log t$ for the variance of $\tilde u(t,\omega,\omega_1)$ under $\tilde \Pm$. The augmented probability space was constructed in such a way that $\tilde u(t,\omega,\omega_1)$ is amenable to Talagrand's inequality. The function $u$ depends on $\omega_j$ for only $O(t^d)$ of the indices $j \in \mathbb{Z}^d$:

\begin{lemma}\label{lem:finitedependenceC} There is a constant $R > 0$ such that 
\[
\tilde u(t,\omega,\omega_1) = \tilde u(t,\phi_j \omega,\omega_1) , \quad \forall\;\; j \in \mathbb{Z}^d, \;\; |j| > R t,\;\; t > 0
\]
holds $\tilde P$ almost surely.
\end{lemma}
{\bf Proof:} Since $|h(\omega_1)| \leq m \sqrt{d} $, this is a consequence of Lemma \ref{lem:finspeed2C}: no approximate optimizer passes through cube $j$, if $R$ is sufficiently large and $|j| > Rt$. \hfill $\Box$

\vspace{0.2in}

In view of Lemma \ref{lem:finitedependenceC}, we may regard $\tilde u$ as a function of no more than $Ct^d + dm^2$ random variables taking values in the set $\{a, b\}$. In this way, Talagrand's inequality (Theorem \ref{talagrands_inequality}) implies that there is a constant $C > 0$, independent of $t > 0$, such that 
\begin{eqnarray} \label{talagrandAug}
\var(\tilde u)\leq  C  \sum_{j\in B_t} \frac{\|\rho_j \tilde u\|_2^2}{1+\log\frac{\|\rho_j \tilde u\|_2}{\|\rho_j \tilde u\|_1}} +  C \sum_{k \in \Theta_t} \frac{\|\rho_k \tilde u\|_2^2}{1+\log\frac{\|\rho_k \tilde u\|_2}{\|\rho_k \tilde u\|_1}}.
\end{eqnarray}
where $B_t$ is the set $B_t = \{ j \in \mathbb{Z}^d\;|\;\; |j| \leq R t \}$, whose cardinality is bounded by $Ct^d$. The norms $\| \cdot \|_2$ and $\| \cdot \|_1$ refer to the $L^2(\tilde \Omega,\tilde P)$ and $L^1(\tilde \Omega,\tilde P)$ norms, respectively.  Observe that if $k \in \Theta_t$, then $\rho_k \tilde f$ corresponds to translation of the random environment:
\[
\rho_k \tilde u = \frac{\sigma_k \tilde u - \tilde u}{2} = \frac{\tilde u(t,\omega,\phi_k \omega_1) - \tilde u(t,\omega,\omega_1)}{2}.
\]
If $j \in B_t$, then $\rho_j \tilde u$ corresponds to a local change in the random environment over the cube $Q_j$:
\[
\rho_j \tilde u = \frac{\sigma_j \tilde u - \tilde u}{2} = \frac{\tilde u(t,\phi_j \omega, \omega_1) - \tilde u(t,\omega,\omega_1)}{2}.
\]

Let us first consider the second sum in (\ref{talagrandAug}). 
We will show that this sum is $O(t^{2\zeta})$.

\begin{lemma} \label{lem:second_sumC}
There is a constant $C >0 $ such that
\begin{eqnarray}
\sum_{k \in \Theta_t} 
\frac{\|\rho_k \tilde u\|_2^2}{1+\log\frac{\|\rho_k \tilde u\|_2}{\|\rho_k \tilde u\|_1}} 
\leq C t^{2 \zeta} \label{Thetasumbound}
\end{eqnarray}
holds for all $t > 1$.
\end{lemma}
{\bf Proof:} Since there are only $|\Theta_t| = m^2 \leq t^{2\zeta}$ terms in the sum and since
\[
1+\log\frac{\|\rho_k \tilde u\|_2}{\|\rho_k \tilde u\|_1} \geq 1,
\]
the lemma will follow from a uniform bound on $\|\rho_k \tilde u\|_2$. By definition of $h(\omega_1)$, we know that $|h(\phi_k \omega_1) - h(\omega_1)| \leq 1$. So, by Lemma \ref{lem:fftilde}, we have$|\tilde u(t,\omega,\omega_1) - \tilde u(t,\omega,\phi_k \omega_1)| \leq C |h(\phi_k \omega_1) - h(\omega_1)| \leq C$ holds $\tilde P$ almost surely, for all $k \in \Theta_t$, $t \geq 1$.  \hfill \qed

\vspace{0.2in}

Having established (\ref{Thetasumbound}), we now consider the first sum in (\ref{talagrandAug}).

\begin{prop}\label{lem:first_sumC}
 There is a constant $C >0 $ such that
\begin{eqnarray}
\sum_{j\in B_t} \frac{\|\rho_j \tilde u\|_2^2}{1+\log\frac{\|\rho_j \tilde u\|_2}{\|\rho_j \tilde u\|_1}}
\leq C \frac{t}{\log t} \label{B_nsumbound}
\end{eqnarray}
holds for all $t > 1$.
\end{prop}

Since we may have $\alpha \neq \beta$, we will make use of the following fact, proved in the appendix:

\begin{lemma} \label{lem:PAcomp}
Let $C'=\min\left\{\frac{\alpha}{\beta},\frac{\beta}{\alpha}\right\}$ and $C''=\max\left\{\frac{\alpha}{\beta},\frac{\beta}{\alpha}\right\}$. For any measureable set $A \subset \Omega$, 
\begin{equation}
C'\mathbb P(A)\leq \mathbb P(\phi_jA)\leq C''\mathbb P(A) \label{PAcomp}
\end{equation}
holds for all $j \in \mathbb{Z}^d$.  Also, for every nonnegative integrable $\psi$, we have
\begin{eqnarray}\label{instead_of_symmetry}
C'\mathbb E(\psi\circ \phi_j)\leq \mathbb E(\psi)\leq C''\mathbb E(\psi\circ \phi_j).
\end{eqnarray}
\end{lemma}

\noindent{\bf Proof of Proposition \ref{lem:first_sumC}.}  Let us begin by estimating $\|\rho_j \tilde u\|_2^2$. By Lemma \ref{lem:PAcomp}, we have
\[
\|\rho_j \tilde u\|_2^2  =  \tilde \expE\left[(\sigma_j \tilde u-\tilde u)^21_{\sigma_j\tilde u>\tilde u} \right] + \tilde \expE\left[(\sigma_j \tilde u-\tilde u)^21_{\sigma_j\tilde u <\tilde u} \right] \leq C\tilde \expE\left[(\sigma_j \tilde u-\tilde u)^21_{\sigma_j\tilde u< \tilde u} \right].
\]
Recalling the definition (\ref{Ajdef}), let $\tilde {\mathcal{I}}_j \subset \tilde \Omega$ be the event that $Q_j$ is an important cube in the shifted environment:
\[
\tilde {\mathcal{I}}_j  =   \{ (\omega,\omega_1) \in \tilde \Omega \;|\; \tau_{h(\omega_1)} \omega \in {\mathcal{I}}_j \}.
\]
Because of (\ref{PsigA1}), the event $\{ \sigma_j \tilde u < \tilde u\}$ is contained in the event $\tilde {\mathcal{I}}_j$. So, we have 
\be
\|\rho_j \tilde u\|_2^2  \leq  C\tilde \expE\left[(\sigma_j \tilde u-\tilde u)^21_{\tilde {\mathcal{I}}_j} \right]. \label{pj2}
\ee

The difference $|\sigma_j\tilde u - \tilde u|$ could be large in some cases, even on the event $\tilde {\mathcal{I}}_j$, so we will distinguish a few possible scenarios. Let $\tilde {\mathcal{I}}_j^+ \subset \tilde {\mathcal{I}}_j$ denote the event that cube $Q_j$ is {\em very important} in the shifted environment:
\[
\tilde {\mathcal{I}}^+_j = \{ (\omega,\omega_1) \in \tilde {\mathcal{I}}_j \;|\;  \tau_{h(\omega_1)} \omega \in {\mathcal{I}}_j^+ \}.
\]
Similarly, let $\tilde {\mathcal{I}}^-_j = \tilde {\mathcal{I}}_j \setminus \tilde {\mathcal{I}}_j^+$ be the event that the cube $Q_j$ is important but not very important.  Since $\omega \mapsto \tau_{h(\omega_1)} \omega$ is measure preserving on $\Omega$, we have
\[
\tilde \Pm(\{ (\omega,\omega_1) \in \tilde \Omega \;|\; (\tilde u - \sigma_j \tilde u)^2 1_{\tilde {\mathcal{I}}_j^-} > C_0 \}) = \Pm(\{ \omega \in \Omega \;|\; (u - \sigma_j  u)^2 1_{{\mathcal{I}}_j^-} > C_0\}).
\]
Consequently, from Lemma \ref{lem:notveryimportant} and (\ref{pj2}) we have
\br
\|\rho_j \tilde u\|_2^2  & \leq &  C\tilde \expE\left[(\sigma_j \tilde u-\tilde u)^21_{\tilde {\mathcal{I}}_j^+} \right] + C\tilde \expE\left[(\sigma_j \tilde u-\tilde u)^21_{\tilde {\mathcal{I}}_j^-} \right] \no \\
& \leq & C\tilde \expE\left[(\sigma_j \tilde u-\tilde u)^21_{\tilde {\mathcal{I}}_j^+} \right] + C C_0 \tilde \Pm(\tilde {\mathcal{I}}_j). \label{rhoj22}
\er
Whether the event $\tilde {\mathcal{I}}_j^+$ has small probability depends on the function $g(y)$, so we distinguish two cases. Let $\tilde G \subset \tilde \Omega$ denote the even that
\[
|\gamma(t) - \gamma(0)| \geq t^{1/4}, \quad \forall\;\gamma \in M_1(\omega,\omega_1).
\]
On this event, all approximate minimizers must travel a distance at least $O(t^{1/4})$ from their starting point $\gamma(0) = 0$.  According to the following lemma, the probability that minimizers travel that far when a cube $Q_j$ is very important must be small.

\begin{lemma} \label{lem:Gtildebound}
There are constants $\kappa_1, \kappa_2 > 0$ such that
\[
\tilde \Pm(\tilde {\mathcal{I}}_j^+ \cap \tilde G) \leq \kappa_1 e^{-\kappa_2 t^{1/4}}, \quad \forall\;\; t > 0, \quad j \in \Zm^d.
\]
\end{lemma}

Therefore, returning to (\ref{rhoj22}) and using the fact that $|\sigma_j \tilde u - \tilde u|\leq O(t)$ and $\tilde {\mathcal{I}}_j^+ = (\tilde {\mathcal{I}}_j^+ \cap \tilde G) \cup (\tilde {\mathcal{I}}_j^+ \cap \tilde G^C)$, we conclude
\[
\|\rho_j \tilde u\|_2^2  \leq C t^2  e^{-\kappa_2 t^{1/4}} + C \tilde \Pm(\tilde {\mathcal{I}}_j) + C\tilde \expE\left[(\sigma_j \tilde u-\tilde u)^21_{\tilde {\mathcal{I}}_j^+ \cap \tilde G^C}\right].  
\]
Hence, 
\begin{eqnarray}\label{ineq:main_bound}
\sum_{j\in B_t}\|\rho_j \tilde u\|_2^2 &\leq & C |B_t| t^2 e^{- \kappa_2 t} + C \tilde \expE\left[\sum_{j\in B_t} 1_{ \tilde {\mathcal{I}}_j}\right] + C\tilde  \expE\left[\sum_{j\in B_t} (\sigma_j \tilde u-\tilde u)^21_{\tilde {\mathcal{I}}_j^+ \cap \tilde G^C}\right].   
\end{eqnarray}
With probability one, the sum $\sum_{j\in B_t} 1_{\tilde {\mathcal{I}}_j}$ is bounded by $O(t)$ because there can be at most $O(t)$ important cubes, as the total number of cubes visited is $O(t)$, by Lemma \ref{lem:finspeed2C}.

The last term in (\ref{ineq:main_bound}) is bounded as follows. First,
\begin{lemma} \label{lem:confinement}
There are $\kappa_1,\kappa_2 > 0$ such that
\[
\tilde \Pm\left(\{ (\omega,\omega_1) \;|\;\;(\sigma_j \tilde u -  \tilde u)^2 1_{ \tilde G^C} > C_0 t^{1/2} \}\right) \leq \kappa_1 e^{-\kappa_2 t^{1/4}} 
\]
holds for all $t > 1$ and all $j \in \Zm^d$.
\end{lemma}
Furthermore, if $\omega \in  {\mathcal{I}}_j^+$, then $\omega \notin {\mathcal{I}}_k^+$ for any $k \neq j$, since $Q_j$ must be the only important cube. Therefore, since $|\sigma_j \tilde u - \tilde u|\leq |b-a|t$ always holds (by (\ref{sigdif})), we must have 
\[
\tilde \expE\left[\sum_{j\in B_t} (\sigma_j \tilde u-\tilde u)^21_{\tilde {\mathcal{I}}_j^+ \cap \tilde G^C } \right] \leq C_0 t^{1/2} + |b - a|^2 t^2 \kappa_1 e^{-\kappa_2 t^{1/4}}.
\]
Considering (\ref{ineq:main_bound}), we have now shown that there is a constant $C' > 0$ for which
\be
\sum_{j\in B_t}\|\rho_j \tilde u\|_2^2 \leq C't \label{ineq:sum}
\ee
holds for all $t > 1$.

Next we consider the denominator in (\ref{B_nsumbound}). We show that there is a constant $C''>0$ such that
\begin{eqnarray}
 \log\frac{\|\rho_j \tilde u\|_2}{\|\rho_j\tilde u\|_1}&\geq&  C''\log t, \label{ineq:logs}
\end{eqnarray}
for all $t > 1$. By the Cauchy-Schwarz inequality we see that
\[
\norm{\rho_j \tilde u}_1 = \norm{\rho_j \tilde u\cdot 1_{\sigma_j \tilde  u\neq \tilde u}}_1 \leq  \norm{ \rho_j \tilde u }_2\cdot \sqrt{\mathbb P(\sigma_j \tilde u\neq  \tilde u)}.
\]
Since $\sigma_j \sigma_j u = u$, Lemma \ref{lem:PAcomp} implies
\[
\tilde \Pm(\sigma_j \tilde u \neq \tilde u) = \tilde \Pm(\sigma_j \tilde u> \tilde u) + \tilde \Pm(\sigma_j \tilde u< \tilde u) \leq  (1 + C'')\tilde \Pm(\sigma_j \tilde u<\tilde u).
\]
Hence, 
\begin{equation}
\frac{\|\rho_j \tilde u\|_2}{\|\rho_j \tilde u\|_1} \geq  \frac{1}{\sqrt{(1 + C'')\tilde \Pm(\sigma_j \tilde u< \tilde u)}}. \label{L2L1ratio}
\end{equation}
Therefore, to bound $\log(\|\rho_j \tilde u\|_2/\|\rho_j \tilde u\|_1)$ from below, we should find an upper bound for $\tilde \Pm(\sigma_j \tilde u< \tilde u)$. Because of (\ref{PsigA}), we know that
\[
\tilde \Pm(\sigma_j \tilde u< \tilde u) \leq \tilde \Pm(\tilde {\mathcal{I}}_j).
\]

To estimate $\tilde \Pm(\tilde {\mathcal{I}}_j)$ we average in $\omega_1$, as was done in \cite{benjamini_kalai_schramm}:
\begin{eqnarray}
\tilde \Pm(\tilde {\mathcal{I}}_j) = \tilde \Em \left[ 1_{\tilde {\mathcal{I}}_j} \right] & = & \tilde \Em \left[ \tilde \Em\left[ 1_{\tilde {\mathcal{I}}_j} \; | \; \omega \right] \right] \nonumber \\
& = & \Em \left[ \sum_{z \in [0,m-1]^d} \tilde \Em \left[ 1_{\tilde {\mathcal{I}}_j} \;| \;\omega , h(\omega_1) = z \right] P_1(h(\omega_1) = z) \right].
\end{eqnarray}
Observe that $(\omega,\omega_1) \in \tilde {\mathcal{I}}_j$ if and only if there is $\delta > 0$ such that
\[
\pi_{j}(\gamma) > 0 \quad \text{for all} \;\; \gamma \in \tilde M_\delta(\omega,\omega_1),
\]
which holds if and only if for some $\delta > 0$
\[
\pi_{j - z}(\gamma) > 0 \quad \text{for all}\; \; \gamma \in M_\delta(t,\tau_z \omega), \;\;\; z = h(\omega_1).
\]
So, for ${\mathcal{I}}_j$ defined by (\ref{AjdefC}), we have
\begin{eqnarray}
\tilde \Em \left[ 1_{\tilde {\mathcal{I}}_j} \;|\; \omega , h(\omega_1) = z \right] = 1_{{\mathcal{I}}_{j-z}}(\tau_z \omega).
\end{eqnarray}
By Lemma \ref{lemma3_benjamini_kalai_schramm}, we also know that $P_1(h(\omega_1) = z) \leq C m^{-d}$ for any $z \in [0,m-1]^d$. Therefore, 
\begin{eqnarray}
\tilde \Pm(\tilde {\mathcal{I}}_j) & \leq & C m^{-d}  \Em \left[ \sum_{z} \tilde \Em \left[ 1_{{\mathcal{I}}_j} \;|\; \omega , h(\omega_1) = z \right]  \right]  = C m^{-d}  \Em  \left[ \sum_{z} 1_{{\mathcal{I}}_{j-z}}(\tau_z \omega) \right] \nonumber \\
&=& C m^{-d}  \Em  \left[ \sum_{z} 1_{{\mathcal{I}}_{j-z}}(\omega) \right]. \nonumber 
\end{eqnarray}
The last equality follows from the stationarity of $\Pm$ with respect to $\tau_z$. Now, given $\omega \in \Omega$, the sum
\[
\sum_{z \in [0,m-1]^d} 1_{{\mathcal{I}}_{j-z}}(\omega)
\]
counts the number of important cubes within the box $j - [0,m-1]^d$. These cubes are visited by all paths $\gamma \in M_\delta(\omega)$ for some $\delta > 0$ sufficiently small. Hence, $\tilde \Pm(\tilde {\mathcal{I}}_j) \leq C m^{-d} \Em[ \# \Lambda_j]$ where
\[
\Lambda_j = \bigcup_{n \geq 1} \left \{ k \in \mathbb{Z}^d \;|\;\; j - k \in [0,m-1]^d, \;\;\pi_k(\gamma) > 0\;\;\forall \gamma \in M_{1/n}(\omega) \right \}.
\]
We may interpret the random variable $\# \Lambda_j$ as the number of important cubes in the box $j - [0,m-1]^d$. 

Obviously we have the trivial bound $ \# \Lambda_j \leq O(t)$. This is because each path $\gamma \in M_\delta(\omega)$ has length $O(t)$, by Lemma \ref{lem:finspeed2C}, so $\pi_k(\gamma) > 0$ for at most $O(t)$ indices $k$. Therefore,
$$
\mathbb E[\#\Lambda_j]\leq t\leq (m+1)^{1/\zeta}.
$$
If $d \geq 3$, we may choose $\zeta \in (1/d,1/2)$, so that  
\[
\tilde \Pm(\sigma_j \tilde u <  \tilde u) \leq \tilde \Pm(\tilde {\mathcal{I}}_j) \leq Cm^{-d} \Em[ \# \Lambda_j]  \leq C m^{- d + 1/\zeta} \leq C t^{1 - d\zeta}
\]
with $1 - d\zeta < 0$. If $d=2$, we need the following:

\begin{lemma} \label{lem:Lambdabound2C}
Let $d = 2$ and assume the non-degeneracy condition (\ref{nondeg2}) holds. Then for each 
$p\in \left(\frac{\nu-1+\zeta}{\zeta \nu},2\right)$ there exists a constant $C >0$ such that  $ \# \Lambda_j \leq C m^{p}$ holds with probability one, for all $j \in \mathbb{Z}^d$.
\end{lemma}

So, for $d=2$ we still have $\tilde \Pm(\sigma_j \tilde u < \tilde u) \leq  Cm^{-d} \Em[ \# \Lambda_j]  \leq C n^{(p-2)/\zeta}$, with $(p-2)/\zeta < 0$. Therefore, returning to (\ref{L2L1ratio}) we conclude that there is a constant $C > 0$ such that
\begin{eqnarray}\label{ineq:LowerBoundL2L1}
\log \frac{\|\rho_j \tilde u\|_2}{\|\rho_j \tilde u\|_1} \geq  C \log t
\end{eqnarray}
holds for all $t$ sufficiently large. Therefore, the proof of Proposition \ref{lem:first_sumC} is reduced to a proof of Lemma \ref{lem:Gtildebound}, Lemma \ref{lem:confinement} and, in case $d=2$, Lemma \ref{lem:Lambdabound2C}. These are proved in the next section. \qed

\vspace{0.2in}

Theorem \ref{fk_model_variance} now follows immediately from (\ref{talagrandAug}), Lemma \ref{lem:second_sumC}, and Proposition \ref{lem:first_sumC}.

\vspace{0.2in}

\section{Proofs of the technical estimates} \label{sec:tech}

\subsection*{Proof of Lemma \ref{lem:Gtildebound}.}  

Observe that 
\[
\tilde \Pm(\tilde {\mathcal{I}}_j^+ \cap \tilde G) = \Pm({\mathcal{I}}^+_j \cap G)
\]
where $G \subset \Omega$ is the event for which 
\begin{equation} \label{gammamove}
|\gamma(t) - \gamma(0)| \geq  t^{1/4}, \quad \forall \;\;\gamma \in M_\delta(t,\omega).
\end{equation}
We will show that on the event ${\mathcal{I}}^+_j \cap G$, any approximate optimizer $\gamma \in M_\delta(\omega)$ must touch a set of $O(t^{1/4})$ cubes which are almost uniformly spaced on a straight line segment of length $O(t)$ and on each of those cubes we have $V(x,\omega) = b$. Such an event can occur only with small probability.

Suppose $\omega \in {\mathcal{I}}^+_j \cap G$, and let $\gamma \in M_\delta(\omega)$ be such that (\ref{gammamove}) holds. Let $[t_1,t_2] \subset [0,t]$ be the smallest interval containing all $s$ for which $\gamma(s) \in Q_j$. Hence, $\omega(\gamma(s)) = b$ for all $s \notin [t_1,t_2]$, since $\omega \in {\mathcal{I}}_j^+$.  Since $\omega \in G$, we know that $|\gamma(t) - \gamma(0)| \geq  t^{1/4}$, which means that either
\[
|\gamma(t_1) - \gamma(0)| > \frac{t^{1/4}}{3}  \quad \quad \text{or} \quad \quad |\gamma(t) - \gamma(t_2)| > \frac{t^{1/4} }{3} 
\]
must hold, because $\gamma(t_1), \gamma(t_2) \in \overline{Q_j}$.  Let us assume that $|\gamma(t) - \gamma(t_2)| > (t^{1/4})/3$ holds; the other case is treated in a similar manner.

First, since $\omega(\gamma(s)) = b$ for all $s \in (t_2,t]$, we may assume that $\gamma$ is a straight line between $\gamma(t_2)$ and $\gamma(t)$. Specifically, by redefining $\gamma$ slightly, we may assume that 
\[
\gamma(s) = \gamma(t_2) + \frac{\gamma(t) - \gamma(t_2)}{t - t_2} (s - t_2), \quad \forall \;\; s \in [t_2,t],
\]
for otherwise, $\gamma$ would not be an optimal path. This follows from (\ref{lowerLsumC}).

Next, given points $\gamma(t_2)$ and $\gamma(t)$, there is a unique pair $x_{t_2}, x_t \in \mathbb{Z}^d$ such that
\[
\gamma(t_2) = x_{t_2} + y_{t_2}, \quad \quad \gamma(t) = x_t + y_t
\]
for some $y_{t_2}, y_t \in [0,1)^d$.  Therefore, if we define the linear path
\[
\hat \gamma(s) = x_{t_2} + \frac{x_t - x_{t_2}}{t-t_2}(s - t_2), \quad s \in [t_2,t]
\]
we have $|\gamma(s) - \hat \gamma(s)| \leq 2 \sqrt{d}$ for $s \in [t_2,t]$. Therefore, for each $s \in [t_2,t]$ there must be a cube $Q_{\ell}$ such that $\mbox{dist}\,(\hat \gamma(s), Q_\ell) \leq 2 \sqrt{d}$ and $\omega(Q_\ell) = b$.  For $y \in \mathbb{R}^d$, let $B_{y}$ denote the event that there is at least one cube $Q$ such that $\mbox{dist}\,(Q,y) \leq 2 \sqrt{d}$ and $\omega(Q) = b$.  Then $\Pm(B_y) = 1 - \Pm(B_y^C) \leq 1  - \alpha^{C_3} < 1$, for a constant $C_3 > 0$ that depends only on the dimension $d$. Moreover, if $|y - z| > 5 \sqrt{d}$, then $B_y$ and $B_z$ are independent events. Therefore, for fixed times $t_2 < t$ and a fixed pair of points $x_{t_2}, x_t \in \mathbb{Z}^d$ satisfying $|x_{t_2} - x_{t}| \geq  t^{1/4}/2$ we have
\[
\Pm\left( \bigcap_{s \in [t_2,t]} B_{\hat \gamma(s)} \right) \leq (1  - \alpha^{C_3})^{C_4  t^{1/4}}
\]
for some constant $C_4 > 0$ independent of $\epsilon$.

By Lemma \ref{lem:finspeed2C}, we know there is a constant $R>0$ such that $|\gamma(s) - \gamma(0)|\leq t R$ for all $s \in [0,t]$. There are at most $O(t^{2d})$ possible pairs $x_{t_2}, x_t \in \mathbb{Z}^d$ satisfying $|x_{t_2} - \gamma(0)| \leq R t$ and $|x_t - \gamma(0)| \leq R t$ and $|x_{t_2} - x_t| \geq  t^{1/4}$. Therefore, we conclude that
\begin{equation}
\Pm \left(  {\mathcal{I}}_j^+ \cap  G \right) \leq  O(t^{2d}) (1  - \alpha^{C_3})^{C_4  t^{1/4}}.
\end{equation} The last inequality immediately implies the lemma.
\hfill $\Box$ 

\vspace{0.2in}


\subsection*{Proof of Lemma \ref{lem:confinement}} 
Because $\omega \mapsto \tau_{h(\omega_1)} \omega$ is measure preserving on $\Omega$, we have
\[
\tilde \Pm\left(\{ (\omega,\omega_1) \;|\;\;(\sigma_j \tilde u -  \tilde u)^2 1_{ \tilde G^C} > C_0 t^{1/2} \}\right)= \Pm\left(\{ \omega \;|\;\;(\sigma_j u -  u)^2 1_{  G^C} > C_0 t^{1/2} \}\right)
\]
where the event $G \subset \Omega$ is defined by (\ref{gammamove}). So, on the event $G^C$ we know there is $\gamma \in M_\delta(\omega)$ such that
\begin{equation}
|\gamma(t) - \gamma(0)| <  t^{1/4}.
\end{equation}
Let $B_{r}(x)$ denote the ball of radius $r > 0$ centered at $x$. We may assume that there are at least two indices $j,k \in \mathbb{Z}^d \cap B_{t^{1/4}}(0)$ such that $\omega_j = a$ and $\omega_k = a$. This is because the event that $\omega_\ell = a$ for at most one of the cubes contained in $B_{t^{1/4}}(0)$ has probability less than $O(\beta^{N_t})$ where $N_t \geq C t^{1/4}$ is the number of cubes contained in $B_{t^{1/4}}(0)$. 

Let $\gamma \in M_\delta(\omega)$ with $|\gamma(t) - \gamma(0)| \leq t^{1/4}$. Then
\[
u(t,\omega) \leq g(\gamma(t)) - at + \delta.
\]
Suppose $\omega_k = a$ for some $k \neq j$ and $k \in B_{t^{1/4}}(0)$. Let $x_k \in Q_k$, so that $V(x_k,\omega) = a$. Define the path $\hat \gamma$ by
\br
\hat \gamma(s) = \left\{ \begin{array}{cc} \gamma(0) + s\frac{x_k - \gamma(0)}{ t^{1/4}},& s \in [0,t^{1/4}] \\ x_k, & s \in [t^{1/4},t - t^{1/4}] \\
x_k + (s - t + t^{1/4})\frac{\gamma(t) - x_k}{t^{1/4}},& s \in [t-t^{1/4},t]. \end{array} \right. 
\er
Then
\[
\sigma_j u(t,\omega) \geq g(\hat \gamma(t)) - \mathcal{L}(\hat \gamma,\omega) \geq g(\gamma(t)) - a(t - 2t^{1/4} - b2t^{1/4} - 2t^{1/4}\max_{|z| \leq 1} K(z).
\]
Therefore,
\[
u(t,\omega) - \sigma_j u(t,\omega) \leq (a - b)2t^{1/4} + 2 t^{1/4} \max_{|z| \leq 1} K(z).
\]
Hence $(u - \sigma_j u)^2 \leq C_0 t^{1/2}$ except possibly on a set of probability less than $O(\beta^{N_t})$.

\hfill \qed


\subsection*{Proof of Lemma \ref{lem:Lambdabound2C} for $d=2$}
Assuming the non-degeneracy condition (\ref{nondeg2}), we may choose real numbers $\nu>1$ and $\varepsilon_0 > 0$ such that $K(q)\geq |q|^\nu$ for all $q$ that satisfy $|q|<\varepsilon_0$.  Having fixed $\zeta \in (\frac{\nu-1}{2\nu-1}, \frac{1}{2})$, we see that $\frac{\nu-1 + \zeta}{\zeta \nu} < 2$. So, we may choose $p \in (\frac{\nu-1 + \zeta}{\zeta \nu}, 2)$.

Arguing by contradiction, we assume that $\#\Lambda_j > m^p$: there are more than $m^p$ important cubes within the box $B_j=j-[0,m-1]^d$. Fix $\delta > 0$ small. Consider a path $\gamma \in M_\delta(n,\omega)$. Let $[t_1,t_2]$ be the smallest interval containing all $s$ for which $\gamma(s) \in B_j$. Hence $|t_2 - t_1| \geq C m^p$.   Choose any one of the important cubes in $B_j$ and let $x_c$ denote its center point. Let us define a modified path $\hat\gamma$ as follows: 
\begin{enumerate}
\item[(i)] For $s \in [0,t_1] \cup [t_2,t]$, let $\hat \gamma(s)=\gamma(s)$.
\item[(ii)] For $s \in [t_1,t_2]$ let 
\[\hat \gamma(s)=\left\{\begin{array}{ll}\gamma(t_1)+ (s-t_1)\cdot \frac{x_c-\gamma(t_1)}{m}, & s\in [t_1,t_1+m],\\
x_c, & s\in [t_1 + m, t_2-m],\\
x_c+(s-t_2+m)\cdot\frac{\gamma(t_2)-x_c}{m}, & s\in[t_2-m, t_2].
 \end{array}\right.
\]
\end{enumerate}
We have the bound
\begin{eqnarray}
\label{ineq:Contradiction_ready}
\mathcal L(\gamma)-\mathcal L(\hat\gamma)&\geq& -2m\cdot \max_{|q|\leq 1} K(q) - 2m(b-a)+\int_{t_1}^{t_2} K(\gamma'(s)) \,ds.
\end{eqnarray}
We will prove that for sufficiently large $m$ the right side of (\ref{ineq:Contradiction_ready}) is larger than $\delta > 0$, contradicting the fact that $\gamma \in M_\delta(n,\omega)$.

Let us denote by $J_0 \subset [t_1,t_2]$ the set of times for which $
|\gamma'(s)|\geq \varepsilon_0$. Therefore, we may assume
\be
\int_{J_0}K(\gamma'(s))\,ds \leq  2 m \left( \max_{|q|\leq 1}K(q)+(b-a)\right) + \delta, \label{J0Kintegral}
\ee
for otherwise the right side of (\ref{ineq:Contradiction_ready}) would be larger than $\delta$. Let $J_1=[t_1,t_2]\setminus J_0$; for these times $s \in J_1$ we have $|\gamma'(s)| \leq \epsilon_0$.  From (\ref{ineq:Contradiction_ready}) we also obtain:
\begin{eqnarray}
\nonumber \mathcal L(\gamma)-\mathcal L(\hat\gamma) &\geq & -2m\left(\max_{|q|\leq 1}K(q)+b-a \right) +\int_{J_1} K(\gamma'(r))\,dr \no \\
&\geq& -2m\left(\max_{|q|\leq 1}K(q)+b-a \right) +\int_{J_1} |\gamma'(r)|^\nu \,dr \nonumber \\
&\geq& -2m\left(\max_{|q|\leq 1}K(q)+b-a \right) +|J_1|\cdot \left(\frac1{|J_1|}\int_{J_1} |\gamma'(r)| \,dr\right)^\nu . 
\label{ineq:Contradiction_ready2}
\end{eqnarray}
In the last line we applied Jensen's inequality. 
We will now prove that there exists a real number $\epsilon_1>0$ such that 
\begin{eqnarray}\label{ineq:BoundOnLength}
\int_{J_1}|\gamma'(r)|\,dr \geq \epsilon_1 m^p.
\end{eqnarray}
By our assumption, the number of important cubes within $B_j$ is more than $m^p$. Let us now paint all these cubes in $2^d$ colors so that no two cubes share the same color. By the pigeon-hole principle there are at least $m^p 2^{-d}$ important cubes having the same color. The distance between two cubes of the same color is at least $1$, hence we have $\int_{t_1}^{t_2}|\gamma'(r)|\,dr \geq m^p2^{-d}$. Therefore, since there is $C > 0$ such that $|\gamma'(s)| \leq C K(\gamma'(s))$ for all $s \in J_0$, we have
\begin{eqnarray}
m^p2^{-d}&\leq& \int_{J_1}|\gamma'(r)|\,dr+\int_{J_0}|\gamma'(r)|\,dr \no \\
& \leq & \int_{J_1} |\gamma'(r)|\,dr + C \int_{J_0} K(\gamma'(r))\,dr \no \\
& \leq &  \int_{J_1}|\gamma'(r)|\,dr +C (m + \delta).
\end{eqnarray}
In the last step we have applied (\ref{J0Kintegral}). This last inequality implies (\ref{ineq:BoundOnLength}), since $p > 1$.

Now the inequalities (\ref{ineq:Contradiction_ready2}) and (\ref{ineq:BoundOnLength}) imply:
\begin{eqnarray}
\mathcal L(\gamma)-\mathcal L(\hat \gamma)&\geq& 
-2m\left(\max_{|q|\leq 1}K(q)+b-a
\right) +|J_1|\cdot \left(\frac{\epsilon_1 m^p}{|J_1|}\right)^\nu \no \\
&=&-2m\left(\max_{|q|\leq 1}K(q)+b-a
\right) + (\epsilon_1)^{\nu}m^{ p\nu}\cdot|J_1|^{1-\nu} \no \\
&\geq& -2m\left(\max_{|q|\leq 1}K(q)+b-a
\right) + (\epsilon_1)^{\nu}m^{p\nu}\cdot \left(m^{1/\zeta} \right)^{1-\nu}. \label{zetaspdeterm}
\end{eqnarray}
In the last inequality we have used $|J_1| \leq n = m^{1/\zeta}$. If we have
\[
p > \frac{\nu + \zeta - 1}{\zeta \nu},
\]
then $p\nu + (1-\nu)/\zeta > 1$. In this case, the right side of (\ref{zetaspdeterm}) is positive, and larger than $\delta$, for $t$ sufficiently large. Since this contradicts the approximate optimality of $\gamma \in M_\delta$, we must have $\# \Lambda_j \leq m^p$. \hfill \qed

\section{Appendix} \label{sec:append}

\noindent {\bf Proof of Lemma \ref{lem:PAcomp}:} The bounds in (\ref{PAcomp}) follow from the fact that $\Pm$ is the product measure on $\Omega = \{a,b\}^{\mathbb{Z}^d_n}$, with $\Pm(\omega(j) = a) = \alpha$ and $\Pm(\omega(j) = b) = \beta$. 
For every nonnegative integrable $\psi$, (\ref{PAcomp}) implies
\begin{eqnarray*}
\mathbb E(\psi)&=&\int \psi\,d\mathbb P \leq C''\int \psi\, d\mathbb P\circ\phi_j=\int \psi(\phi_j\omega)\, d\mathbb P=\mathbb E(\psi\circ \phi_j).
\end{eqnarray*}
Similarly $\mathbb E(\psi)\geq C'\mathbb E(\psi\circ \phi_j)$ for all such $\psi$. \hfill $\Box$

\vspace{0.2in}

\noindent {\bf Proof of Theorem \ref{talagrands_inequality}:}  Let us define 
\[
\Delta_j f(\omega)=\left\{\begin{array}{ll} \beta(f(\phi_j\omega) - f(\omega)),& \mbox{ if } \omega_j=a\\
\alpha(f(\phi_j\omega) - f(\omega)),&\mbox{ if }\omega_j=b.\end{array}\right. 
\]
Then Theorem \ref{talagrands_inequality} is a slight modification of the following

\begin{theorem}[\cite{talagrand_94}, Theorem 1.5] \label{talagrands_inequality_original}
There is a constant $C > 0$, such that 
\begin{equation}
\var(f)\leq 
C\cdot \sum_{j\in J}\frac{\|\Delta_j f\|_2^2}{1+\log\frac{\|\Delta_jf\|_2}{\|\Delta_jf\|_1}}. \label{talagrand_ineq0}
\end{equation}
holds for all $f \in L^2(\Omega_J)$.
\end{theorem}

To derive Theorem \ref{talagrands_inequality} from this, we start with elementary observation $$C'|\rho_j f(\omega)|\leq |\Delta_j f(\omega)|\leq C''|\rho_j f(\omega)|$$ for $C'=\min\{2\alpha,2\beta\}$ and $C''=\max\{2\alpha,2\beta\}$. Let $\kappa = \log(C''/C')  \geq 0$. If $\log\frac{\|\rho_jf\|_2}{\|\rho_jf\|_1} \geq 2\kappa$, then 
\[
\log\frac{\|\Delta_j f\|_2}{\|\Delta_jf\|_1} \geq \log\frac{C'}{C''}+\log\frac{\|\rho_jf\|_2}{\|\rho_jf\|_1} \geq \frac{1}{2} \log\frac{\|\rho_jf\|_2}{\|\rho_jf\|_1}.
\] 
Consequently, Theorem \ref{talagrands_inequality_original} implies
\[
\var(f)\leq 
C\cdot \sum_{j\in J}\frac{\|\Delta_j f\|_2^2}{1+\log\frac{\|\Delta_jf\|_2}{\|\Delta_jf\|_1}} \leq  2(C'')^2 C\sum_{j\in J}\frac{\|\rho_jf\|_2^2}{1+\log\frac{\|\rho_jf\|_2}{\|\rho_jf\|_1}}.
\]
On the other hand, if $\log\frac{\|\rho_jf\|_2}{\|\rho_jf\|_1} \in [0, 2\kappa)$, then Theorem \ref{talagrands_inequality_original} implies
\[
\var(f)\leq 
C\cdot \sum_{j\in J}\|\Delta_j f\|_2^2 \leq  (1 + 2\kappa)2(C'')^2 C\sum_{j\in J}\frac{\|\rho_jf\|_2^2}{1+\log\frac{\|\rho_jf\|_2}{\|\rho_jf\|_1}}.
\]
\hfill $\Box$


\end{document}